\documentclass[12pt]{amsart}
\usepackage{amsmath, amssymb, latexsym, amsthm, mathrsfs, color, mathtools, tikz,pgfplots,tikz-cd, graphicx, stackengine, url}
\usepackage[top=2.5cm, bottom=2.5cm, left=2.5cm, right=2.5cm]{geometry}

\stackMath
\newcommand\tsup[2][2]{%
 \def\useanchorwidth{T}%
  \ifnum#1>1%
    \stackon[-1pt]{\tsup[\numexpr#1-1\relax]{#2}}{\hspace{1pt}\scriptstyle\sim}%
  \else%
    \stackon[.5pt]{#2}{\hspace{1pt}\scriptstyle\sim}%
  \fi%
}

\newcommand{\nc}{\newcommand}

\DeclareMathOperator{\Int}{Int}

\newcommand{\sgg}{{\mathsf{S}_1(\Ga,\Ga)}}

\nc{\mc}{\mathcal}
\nc{\thusfar}{\my{--- Edited thus far ---}}
\nc{\lei}{\le^\oo}
\nc{\sqsubs}{\sqsubseteq^*}
\nc{\card}[1]{\left|#1\right|}
\nc{\medcard}[1]{\biggl|\,#1\,\biggr|}
\nc{\smallcard}[1]{|\,#1\,|}
\nc{\bds}{bidirectional $\roth$-scale}
\nc{\bfP}{\mathbf{P}}
\nc{\bfQ}{\mathbf{Q}}
\nc{\bbT}{\mathbb{T}}
\nc{\bbZ}{\mathbb{Z}}
\nc{\bbN}{\mathbb{N}}
\nc{\bbC}{\mathbb{C}}
\nc{\beq}{\begin{equation}}\nc{\eeq}{\end{equation}}
\nc{\mbq}{\mb{?}}
\nc{\mb}[1]{{\mbox{\textbf{#1}}}}
\nc{\nop}{$\times$}
\nc{\fbn}{\!\!\fbox{\!\nop\!}\!\!}
\nc{\yup}{\checkmark}
\nc{\forces}{\Vdash}
\nc{\name}[1]{\dot{#1}}
\nc{\tf}{\my{FINISHED THUS FAR}}
\nc{\FU}{Fr\'echet--Urysohn}
\nc{\gs}{$\gamma$~space}
\nc{\Ga}{\Gamma}
\nc{\Gab}{\Ga_\mathrm{Bor}}
\nc{\Om}{\Omega}

\nc{\smallbinom}[2]{\begin{psmallmatrix} #1\\ #2 \end{psmallmatrix}}
\nc{\bgamma}{\smallbinom{\Om}{\Ga}}
\newcommand{\two}{\{0,1\}}
\nc{\productive}[2]{(#1,\allowbreak #2)^\x}
\nc{\prdct}[1]{#1^\x}
\nc{\Sel}{\mathsf{S}}
\nc{\sset}[2]{\{\,#1 : #2\,\}}
\nc{\smb}[1]{{\!\!\mb{#1}\!\!}}
\nc{\medset}[2]{{\biggl\{\,#1 : #2\,\biggr\}}}
\nc{\smallmedset}[2]{{\bigl\{\,#1 : #2\,\bigr\}}}
\nc{\set}[2]{{\left\{\,#1 : #2\,\right\}}}
\nc{\eseq}[1]{#1_1, \allowbreak #1_2, \allowbreak\dotsc} 
\nc{\seleseq}[1]{#1_1\in \mathcal{#1}_1, \allowbreak #1_2\in \mathcal{#1}_2, \allowbreak\dotsc}
\nc{\cube}{(\Cantor)^\bbN}
\nc{\Match}{\op{Match}}
\nc{\concat}[1]{\hat{\phantom{a}}\langle #1\rangle}
\nc{\poset}{\mathbb{P}}
\nc{\fn}[1]{{\op{Fn}(#1\times\w,2)}}
\nc{\linadd}{\op{linadd}}
\nc{\nonprod}{\non^\x}
\nc{\alephes}{{\aleph_0}}
\nc{\my}[1]{\marginpar{\textcolor{red}{***}}\textcolor{red}{#1}}
\nc{\myb}[1]{\marginpar{\textcolor{blue}{***}}\textcolor{blue}{#1}}
\nc{\later}[1]{{\color{green} #1}}
\nc{\BTs}[1]{{\color{green} #1 (BT)}}
\nc{\Cp}{\op{C}_\mathrm{p}}
\nc{\Bp}{\op{B}_p}
\nc{\Pa}[8]{\bibitem{#1} {#2}, \emph{#3}, {#4} \textbf{#5} ({#6}), {#7}--{#8}.}
\nc{\tPa}[5]{\bibitem{#1} {#2}, \emph{#3}, {#4}, to appear.}
\nc{\sPa}[4]{\bibitem{#1} {#2}, \emph{#3}, {#4}, submitted.}
\nc{\Bc}[9]{\bibitem{#1} {#2}, \emph{#3}, in: \textbf{#4} (#5), #6 #7, #8--#9.}
\nc{\fD}{\mathfrak{D}}
\nc{\fX}{\mathfrak{X}}
\nc{\Onbd}{\Op_{\mathrm{nbd}}} 
\nc{\Omnb}{\Om_{\mathrm{nbd}}} 
\nc{\od}{\mathfrak{od}}
\nc{\Setting}[7]{\xymatrix@R=4pt@C=7pt{#1\ar@{-}[r]&#2\ar@{-}[r]&#3\\&#4\ar@{-}[u]\\
#5\ar@{-}[uu]\ar@{-}[r] & #6\ar@{-}[u]\ar@{-}[r] & #7\ar@{-}[uu]}}
\nc{\mx}[1]{\begin{matrix}#1\end{matrix}}
\nc{\plim}{p\txt{-}\lim}
\nc{\Bgp}{{\Z^\bbN}}
\nc{\Cgp}{{{\Z_2}^\bbN}}
\nc{\Cite}[1]{\textbf{[#1]}}
\nc{\Next}[1]{{#1^+}}
\nc{\cFin}{\mathrm{cF}}
\nc{\scsp}{\text{-scale space}}
\nc{\cfn}{\text{cofinal}\ }
\nc{\Con}{\text{Concentrated}}
\nc{\Lind}{\text{Lindel\"of}\,}
\nc{\con}{\text{-Concentrated}}
\nc{\lind}{\text{-Lindel\"of}\,}
\nc{\ctbl}{\text{countably }\allowbreak}
\nc{\Hur}{\text{Hurewicz}}
\nc{\intvl}[2]{{[#1(#2),\allowbreak #1(#2\!+\!1))}}
\nc{\Bdd}{\mathbf{B}}
\nc{\Dfin}{\mathfrak{D}_\mathrm{fin}}
\nc{\grbl}{{\mbox{\textit{\tiny gp}}}}
\nc{\bbP}{\mathbb{P}}
\nc{\BOfat}{\B_{\Om_{\mathrm{fat}}}}
\nc{\Bgood}{\B_{\mathrm{good}}}
\nc{\compactN}{\cl{\mathbb{N}}}
\nc{\blocks}[2]{\op{cl}_{#2}(#1)}
\nc{\blocksplus}[2]{\op{cl}^+_{#2}(#1)}
\nc{\arx}[1]{\texttt{http://arxiv.org/math/#1}}
\nc{\bq}{\begin{quote}}
\nc{\eq}{\end{quote}}
\nc{\cl}[1]{\overline{#1}}
\nc{\Cl}[2]{\mathrm{cl}_{#1}(#2)}
\nc{\CH}{the Continuum Hypothesis}
\nc{\MA}{Martin's Axiom}
\nc{\Bfat}{\B_\mathrm{fat}}
\nc{\inv}{^{-1}}
\nc{\Cantor}{{\two^\bbN}}
\nc{\bP}{\mathbf{P}}
\nc{\bof}{\op{\fb}}
\nc{\dof}{\op{\fd}}
\nc{\bofF}{\bof(\cF)}
\nc{\sr}[3]{\underset{\mbox{#3}}{\mbox{#1}}}
\nc{\gp}{\binom{\Om}{\Ga}}
\nc{\gpsmall}{\mbox{$\gp$}}
\nc{\gig}{\gimel}
\nc{\gns}{\sone(\Om,\gig)}
\nc{\nsr}[2]{#1}
\nc{\Srg}{{\mathbb{S}}}
\nc{\Srgs}{{\mathbb{S}^*}}
\nc{\NN}{{\bbN^{\bbN}}}
\nc{\ZN}{{\Z^{\bbN}}}
\nc{\NNup}{{\bbN^{\uparrow\bbN}}}
\nc{\NNupb}{{b^{\uparrow\bbN}}}
\nc{\Pof}{\op{P}}
\nc{\PN}{{\Pof(\bbN)}}
\nc{\rothel}[1]{{[#1]^{\mbox{\tiny $\infty$}}}}
\nc{\roth}{{[\bbN]^{\mbox{\tiny $\infty$}}}} 
\nc{\roths}{{[b]^{\mbox{\tiny $\infty$}}}} 
\nc{\Fin}{\mathrm{Fin}}
\nc{\ici}{[\bbN]^{ \infty, \infty}}
\nc{\Inc}{{\compactN^{\uparrow\bbN}}}
\nc{\powInc}[1]{{\big(\Inc\big)^{#1}}}
\nc{\powFin}[1]{{\big(\Fin\big)^{#1}}}
\nc{\powPN}[1]{{\big(\PN\big)^{#1}}}
\nc{\NcompactN}{{\compactN^\bbN}}
\nc{\seq}[2]{\la #1\ra_{#2\in\bbN}}
\nc{\Uarrow}{\smash{\big\uparrow}}
\nc{\LE}{\preccurlyeq}
\nc{\GE}{\succcurlyeq}
\nc{\op}{\operatorname}
\nc{\im}{\op{im}}
\nc{\Span}{\op{span}}
\nc{\maxfin}{\op{maxfin}}
\nc{\ran}{\op{range}}
\nc{\iso}{\cong}
\nc{\Madd}{{\M}^\star}
\nc{\cI}{\mathcal{I}}
\nc{\cJ}{\mathcal{J}}
\nc{\scrA}{\mathscr{A}}
\nc{\scrB}{\mathscr{B}}
\nc{\scrC}{\mathscr{C}}
\nc{\scrD}{\mathscr{D}}
\nc{\scrF}{\mathscr{F}}
\nc{\scrK}{\mathscr{K}}
\nc{\A}{\D\forall}
\nc{\B}{\mathrm{B}}
\nc{\cB}{\mathcal{B}}
\nc{\cZ}{\mathcal{Z}}
\nc{\bB}{\mathbf{B}}
\nc{\BS}{\mathbf{B}(\mathcal{S})}
\nc{\BF}{\mathbf{B}(\mathcal{F})}
\nc{\BU}{\mathbf{B}(\mathcal{U})}
\nc{\cSp}{\mathcal{S}^+}
\nc{\cFp}{\mathcal{F}^+}
\nc{\cUp}{\mathcal{U}^+}

\nc{\BG}{\B_\Ga}
\nc{\BL}{\B_\Lambda}
\nc{\BT}{\B_\Tau}
\nc{\BTstar}{\B_{\Tau^*}}
\nc{\BO}{\B_\Om}
\nc{\DO}{\cD_\Om}
\nc{\KO}{\cK_\Om}
\nc{\CG}{C_\Ga}
\nc{\CL}{C_\Lambda}
\nc{\CT}{C_\Tau}
\nc{\CTstar}{C_{\Tau^*}}
\nc{\CO}{C_\Om}
\nc{\COgp}{C_{\Om^{\grbl}}}
\nc{\CLgp}{C_{\Lambda^{\grbl}}}
\nc{\BOgp}{\B_{\Om}^{\grbl}}
\nc{\BLgp}{\B_{\Lambda^{\grbl}}}
\nc{\sfC}{\mathsf{C}}
\nc{\sfD}{\mathsf{D}}
\nc{\bD}{\mathbf{D}}
\nc{\Tau}{\mathrm{T}}
\nc{\cA}{\mathcal{A}}
\nc{\cK}{\mathcal{K}}
\nc{\cD}{\mathcal{D}}
\nc{\cF}{\mathcal{F}}
\nc{\cS}{\mathcal{S}}
\nc{\cT}{\mathcal{T}}
\nc{\cG}{\mathcal{G}}
\nc{\cY}{\mathcal{Y}}
\nc{\J}{\mathcal{J}}
\nc{\cL}{\mathcal{L}}
\nc{\cM}{\mathcal{M}}
\nc{\cN}{\mathcal{N}}
\nc{\cH}{\mathcal{H}}
\nc{\cO}{\mathcal{O}}
\nc{\Op}{\mathrm{O}}
\nc{\rmA}{\mathrm{A}}
\nc{\rmF}{\mathrm{F}}
\nc{\rmB}{\mathrm{B}}
\nc{\rmD}{\mathrm{D}}
\nc{\rmP}{\mathrm{P}}
\nc{\cC}{\mathcal{C}}
\nc{\cP}{\mathcal{P}}
\nc{\bbQ}{\mathbb{Q}}
\nc{\bbR}{\mathbb{R}}
\nc{\cU}{\mathcal{U}}
\nc{\cQ}{\mathcal{Q}}
\nc{\Un}{\bigcup}
\nc{\cV}{\mathcal{V}}
\nc{\cW}{\mathcal{W}}
\nc{\Z}{{\mathbb Z}}
\nc{\Impl}{\Rightarrow}
\long\def\forget#1\forgotten{\marginpar{\textcolor{green}{Forgetting...}}}
\nc{\ft}{\mathfrak{t}}
\nc{\fb}{\mathfrak{b}}
\nc{\fc}{\mathfrak{c}}
\nc{\fd}{\mathfrak{d}}
\nc{\fg}{\mathfrak{g}}
\nc{\oo}{\infty}
\nc{\fr}{\mathfrak{r}}
\nc{\fk}{\mathfrak{k}}
\nc{\bidi}{\mathfrak{bidi}}
\nc{\fu}{\mathfrak{u}}
\nc{\fh}{\mathfrak{h}}
\nc{\fp}{\mathfrak{p}}
\nc{\fj}{\mathfrak{j}}
\nc{\fs}{\mathfrak{s}}
\nc{\w}{\omega}
\nc{\x}{\times}
\nc{\Iff}{\Leftrightarrow}

\nc{\nin}{\notin}
\nc{\cat}{\hat{\ }}
\nc{\sub}{\subseteq}
\nc{\spst}{\supseteq}
\nc{\sm}{\setminus}
\nc{\as}{\subseteq^*}
\nc{\les}{\le^*}
\nc{\leinf}{\le^{\infty}}
\nc{\leS}{\le_S}
\nc{\leF}{\le_F}
\nc{\leU}{\le_U}
\nc{\rest}{\restriction}

\nc{\la}{\langle}
\nc{\ra}{\rangle}
\nc{\E}{\exists}
\nc{\dom}{\op{dom}}
\nc{\cov}{\op{cov}}
\nc{\add}{\op{add}}
\nc{\addmen}{\add(\Men{})}
\nc{\cof}{\op{cof}}
\nc{\cf}{\op{cf}}
\nc{\non}{\op{non}}
\nc{\unif}{\op{non}}
\nc{\COV}{\op{COV}}
\nc{\ADD}{\op{ADD}}
\nc{\COF}{\op{COF}}
\nc{\NON}{\op{NON}}
\nc{\impl}{\to}
\nc{\Lp}{\mathcal{L_\p}}
\nc{\Wlog}{without loss of generality}
\newtheorem{thm}{Theorem}[section]
\nc{\bthm}{\begin{thm}} \nc{\ethm}{\end{thm}}
\newtheorem{prop}[thm]{Proposition}
\nc{\bprp}{\begin{prop}} \nc{\eprp}{\end{prop}}
\newtheorem{fact}[thm]{Fact}
\nc{\bfct}{\begin{fact}} \nc{\efct}{\end{fact}}
\newtheorem{prob}[thm]{Problem}
\nc{\bprb}{\begin{prob}} \nc{\eprb}{\end{prob}}
\newtheorem{lem}[thm]{Lemma}
\nc{\blem}{\begin{lem}} \nc{\elem}{\end{lem}}
\newtheorem{app}[thm]{Application}
\nc{\bapp}{\begin{app}} \nc{\eapp}{\end{app}}
\newtheorem{claim}[thm]{Claim}
\nc{\bclm}{\begin{claim}} \nc{\eclm}{\end{claim}}
\newtheorem{cor}[thm]{Corollary}
\nc{\bcor}{\begin{cor}} \nc{\ecor}{\end{cor}}
\newtheorem{conj}[thm]{Conjecture}
\nc{\bcnj}{\begin{conj}} \nc{\ecnj}{\end{conj}}
\theoremstyle{definition}
\newtheorem{defn}[thm]{Definition}
\nc{\bdfn}{\begin{defn}} \nc{\edfn}{\end{defn}}
\newtheorem{obs}[thm]{Observation}
\nc{\bobs}{\begin{obs}} \nc{\eobs}{\end{obs}}
\theoremstyle{remark}
\newtheorem{rem}[thm]{Remark}
\nc{\brem}{\begin{rem}} \nc{\erem}{\end{rem}}
\newtheorem{cnv}[thm]{Convention}
\nc{\bcnv}{\begin{cnv}} \nc{\ecnv}{\end{cnv}}
\newtheorem{exam}[thm]{Example}
\nc{\bexm}{\begin{exam}} \nc{\eexm}{\end{exam}}
\nc{\bpf}{\begin{proof}} \nc{\epf}{\end{proof}}
\nc{\be}{\begin{enumerate}}
\nc{\ee}{\end{enumerate}}
\nc{\bi}{\begin{itemize}}
\nc{\bimy}{\my{\begin{itemize}}
\nc{\eimy}{\end{itemize}}}
\nc{\itm}{\item}
\nc{\ei}{\end{itemize}}
\nc{\Subsection}[1]{\goodbreak\subsection*{#1}}

\nc{\sone}{\mathsf{S}_1}
\nc{\sfin}{\mathsf{S}_\mathrm{fin}}
\nc{\ufin}{\mathsf{U}_\mathrm{fin}}
\nc{\Split}{\mathsf{Split}}

\nc{\gone}{\mathsf{G}_1}    \nc{\gfin}{\mathsf{G}_\mathrm{fin}}

\nc{\men}{\sfin(\Op,\Op)}
\nc{\sch}{\ufin(\Op,\Omega)}
\nc{\rothb}{\sone(\Op,\Op)}
\nc{\pmen}{\sfin(\Omega,\Omega)}
\nc{\Rothb}{\sone(\Op,\Op)}

\nc{\prothb}{\sone(\Omega,\Omega)}
\nc{\tU}{{\tilde{U}}}
\nc{\tF}{{\tilde{F}}}
\nc{\tY}{{\tilde{Y}}}
\nc{\tX}{S}
\nc{\dtX}{X\sm S}
\nc{\dt}[1]{{\tsup[2]{#1}}}
\nc{\td}{{\tilde{d}}}
\nc{\tb}{{\tilde{b}}}
\nc{\tz}{{\tilde{z}}}
\nc{\cfd}{\cf(\fd)}
\nc{\msep}{\sfin(\cD,\cD)}
\nc{\rsep}{\sone(\cD,\cD)}
\nc{\cft}{\sfin(\Omega_{\mathbf{0}},\Omega_{\mathbf{0}})}
\nc{\scft}{\sone(\Omega_{\mathbf{0}},\Omega_{\mathbf{0}})}

\nc{\Umen}{U\text{-Menger}}
\nc{\hur}{\ufin(\cO,\Gamma)}
\nc{\tUmen}{\tU\text{-Menger}}

\nc{\Men}{\text{Menger}}
\nc{\Sch}{\text{Scheepers}}

\nc{\aspst}{\prescript{*}{}{\spst}\ }
\nc{\eqs}{=^*}
\nc{\ctblOm}{\Omega_{\mathrm{ctbl}}}
\nc{\GNga}{{\smallbinom{\Om}{\Ga}}}
\nc{\ctblga}{\smallbinom{\ctblOm}{\Ga}}
\nc{\Gaclp}{\Ga_{\mathrm{clp}}}
\nc{\sclpgg}{\sone(\Gaclp,\Gaclp)}

\nc{\sep}{
\vspace{2cm}
\noindent
\begin{minipage}{\textwidth}
	\textcolor{red}{\rule{\textwidth}{1pt}}
\end{minipage}
}

\nc{\bfzero}{\mathbf{0}}

\nc{\Mich}[1]{{#1}_\textrm{M}}

\usepackage[all]{xy}

\title{Unbounded towers and the Michael line topology}

\author[W.~Przybylska]{Wanda Przybylska}
\address{Wanda Przybylska, Institute of Mathematics, Faculty of Mathematics and Natural Science College of Sciences, Cardinal Stefan Wyszy\'nski University in Warsaw, W\'oycickiego 1$\slash$3, 01--938 Warsaw, Poland
}
\email{przybylska.wanda@gmail.com}

\keywords{unbounded tower, $\sgg$, Gerlits--Nagy, $\gamma$-property, $\gamma$-set, selection principles, products.}

\subjclass[2010]{26A03, 
54D20, 03E75, 03E17.
}

\begin{document}

\maketitle

\begin{abstract}
A topological space satisfies $\GNga$ (also known as Gerlits--Nagy's property $\gamma$) if every open cover of the space such that each finite subset of the space is contained in a member of the cover, contains a point-cofinite cover of the space.
A topological space satisfies $\ctblga$ if in the above definition we consider countable covers.
We prove that subspaces of the Michael line with a special combinatorial structure have the property $\ctblga$.
Then we apply this result to products of sets of reals with the property $\GNga$.
The main method used in the paper is coherent omission of intervals invented by Tsaban.
\end{abstract}

\section{Introduction}

By \emph{space} we mean a topological space.
A \emph{cover} of a space is a family of proper subsets of the space whose union is the entire space.
An \emph{open} cover is a cover whose members are open subsets of the space.
An \emph{$\omega$-cover }is an open cover such that each finite subset of the space is contained in a set from the cover.
A \emph{$\gamma$-cover} is an infinite open cover such that each point of the space belongs to all but finitely many sets from the cover.
Given a space, let $\Omega$, $\ctblOm$, $\Gamma$ be the families of $\omega$-covers, countable $\omega$-covers and $\gamma$-covers, respectively.
For families $\cA$ and $\cB$ of covers of a space, the property that every cover in the family $\cA$ has a subcover in the family $\cB$ is denoted by $\binom{ \ \cA \ }{\cB}$.
The property $\mathsf{S}_1(\cA, \cB)$ means that for each sequence $\cU_1, \cU_2,\dotsc \in\cA$ there are sets $U_1\in\cU_1, U_2\in\cU_2,\dotsc$ such that $\sset{U_n}{n\in\bbN}\in\cB$.

Let $\roth$ be the set of infinite subsets of $\bbN$ and $\Fin$ be the set of finite subsets.
For sets $a, b\in\roth$ we say that	$a$ is \textit{almost subset} of $b$, denoted $a\sub^*b$, it the set $a\sm b$ is finite.
A \textit{pseudointersection} of a family of infinite sets is an infinite sets $a$ with $a\sub^*b$ for all sets $b$ in the family.
A family of infinite sets is \textit{centered }if the finite intersections of its elements, are infinite.
Let $\fp$ be the minimal cardinality of a subfamily of $\roth$ that is centered and has no pseudointersecion.

\bdfn[{\cite[Definition~2.2]{unbddtower}}]
Let $\kappa$ be an uncountable ordinal number. A set $X\sub \roth$ with 
$|X|\geqslant \kappa$  is a $\kappa$-\emph{\emph{generalized tower}} if for each function
$a\in\roth$, there are sets $b\in\roth$ and $S\sub X$ with $|S|<\kappa$ such that
$$ x\cap \Un_{n\in b} [a(n), a(n+1))\in\Fin$$
for all sets $x\in X\sm S$.
\edfn

Let $\kappa$ be an uncountable ordinal number.
A set $X\cup\Fin$ is \emph{$\kappa$-generalized tower set} if the set $X$ is a $\kappa$-generalized tower.

The \emph{Michael line} is the set $\PN$, with the topology where the points of the set $\roth$ are isolated, and the neighborhoods of the points of the set $\Fin$ are those induced by the Cantor space topology on $\PN$.

\blem[{\cite[Lemma~1.2.]{miller}}]\label{Fin}
Let $\cU$ be a family of open sets in $\PN$ such that $\cU\in\Omega(\Fin)$.
There are a function $a\in\roth$ and sets $U_1, U_2, \dotsc\in\cU$ such that for each set $x\in \roth$ and all natural numbers $n$:
\[
\text{If } x\cap [a(n), a(n+1))=\emptyset, \text{ then } x\in U_n.
\]
\elem

For a set $U\sub\PN$, let $\Int (U)$ be the interior of the set $U$ in the Cantor space topology on $\PN$.
If $\cU\in\Omega(\Fin)$ is a family of open sets in $\PN$ with the Michael line topology, then $\sset{\Int (U)}{U\in\cU}\in\Omega(\Fin)$. Thus Lemma~\ref{Fin}  holds for a family of open sets with the Michael line topology.

\section{Main result}

For functions $f,g\in\NN$ let $(f\circ g)\in\NN$ be a function such that $(f\circ g)(n):=f(g(n))$ for all natural numbers $n$.

\bthm\label{thm1}
Let $X\cup\Fin$ be a $\fp$-generalized tower set with the Michael line topology.
The space $X\cup\Fin$ satisfies $\ctblga$.
\ethm

\bpf
Let $\cU\in\ctblOm (X\cup\Fin)$ be a family of open sets in $\PN$ with the Michael line topology.
Let $S_1:=\Fin$.
Fix a natural number $k>1$, and assume that the set $S_{k-1}\sub X\cup\Fin$ with $\Fin\sub\S_{k-1}$ and $|S_{k-1}|<\fp$ has been already defined.
Since $|S_{k-1}|<\fp$, there is $\cV\sub\cU$ such that $\cV\in\Gamma(S_{k-1})$.
By Lemma~\ref{Fin}, there are a function $a_k\in[\bbN]^{\infty}$ and sets $U_1^{(k)}, U_2^{(k)},\dotsc\in\cV$ such that for each set $x\in\roth$ and all natural numbers $n$:
\beq \text{If }x\cap [a_k(n), a_k(n+1))=\emptyset, \text{ then } x\in U_n^{(k)}. \tag{$~\ref{thm1}$.$1$} \eeq
Since the set $X$ is $\fp$-generalized tower, there are a set $b_k\in\roth$ and a set 
$S_k\sub X\cup\Fin$ with $S_{k-1}\sub S_k$ and $|S_k|<\fp$ such that
$$ x\cap \bigcup_{n\in b_k} [a_k(n), a_k(n+1))\in\Fin$$
for all sets $x\in X\sm S_k$.
Then $$ \sset{U_{b_k(j)}^{(k)}}{j\in\bbN}\in\Gamma((X\sm S_k)\cup S_{k-1}) .$$ 
There is a function $a\in\roth$ such that for each natural number $k$, we have
$$ |(a_k\circ b_k)\cap [a(n), a(n+1))|\geqslant 2,$$
for all but finitely many natural numbers $n$.
Since the set $X$ is $\fp$-generalized tower, there are a set $b\in\roth$ and a set 
$S\sub X$ with  $|S|<\fp$ such that
\beq x\cap \bigcup_{n\in b} [a(n), a(n+1))\in\Fin  \tag{$~\ref{thm1}$.$2$} \eeq
for all sets $x\in X\sm S$.
We may assume that $\bigcup_{k} S_k\sub S$.
The sets 
\beq  c_k:=\sset{i\in b_k}{[a_k(i), a_k(i+1))\sub \bigcup_{n\in b}[a(n), a(n+1))} \tag{$~\ref{thm1}.3$} \eeq
are infinite for all natural numbers $k$.
Thus, 
$$ \sset{U^{(k)}_{c_k(j)}}{j\in\bbN} \in \Gamma((X\sm S_k)\cup S_{k-1}).$$
Since the sequence of the sets $S_k$ is increasing, we have $X = \Un_k (X\sm S_k)\cup S_{k-1}$
and each point of $X$ belongs to all but finitely many sets $(X\sm S_k)\cup S_{k-1}$.
For each point $x\in S$, define
 $$ g_x(k):= \begin{cases}
0 & x\notin (X\sm S_k)\cup S_{k-1},\\
\min\sset{j}{x\in\bigcap_{i\geqslant j} U^{(k)}_{c_k(i)}} & x\in (X\sm S_k)\cup S_{k-1}.
\end{cases}$$
Since $|S|<\fp$, there is a function $g\in\bbN^{\bbN}$ with $\sset{g_x}{x\in S}\leqslant^* g$ and
 \beq  a_k(c_k(g(k)+1))<a_{k+1}(c_{k+1}(g(k+1)))  \tag{$~\ref{thm1}.4$} \eeq
for all natural numbers $k$.
Let $$ \cW_k:=\sset{U^{(k)}_{c_k(j)}}{j\geqslant g(k)}$$
for all natural numbers $k$.
Then $\cW_1, \cW_2, \dotsc\in\Gamma(S)$.
We may assume that families $\cW_k$ are pairwise disjoint.
Since properties $\ctblga$ and $\mathsf{S}_1(\ctblOm, \Gamma)$ are equivalent, the set $S$ satisfies $\mathsf{S}_1(\ctblOm, \Gamma)$.
Then there is a function $h\in\bbN^{\bbN}$ such that $g\leqslant h$ and 
 $$ \sset{U^{(k)}_{c_k(h(k))}}{k\in\bbN}\in\Gamma(S).$$
Fix a set $x\in X\sm S$.
By ($~\ref{thm1}.3$), for each natural number $k$, we have
 $$ \Un_{n\in c_k}[a_k(n), a_k(n+1))\sub \Un_{n\in b} [a(n), a(n+1)).$$ 
By ($~\ref{thm1}.2$), ($~\ref{thm1}.4$) and the fact that $g\leqslant h$, the set $x$ omits all but finitely many intervals
 $$ [a_k(c_k(h(k))), a_k(c_k(h(k))+1) ).$$
 By ($~\ref{thm1}.1$), we have
  $$ \sset{U^{(k)}_{c_k(h(k))}}{k\in\bbN}\in\Gamma(X\sm S)).$$
  Then
   $$ \sset{U^{(k)}_{c_k(h(k))}}{k\in\bbN}\in\Gamma(X\cup\Fin).$$
\epf

\section{Applications}

For spaces $X$ and $Y$, let $X\sqcup Y$ be the \emph{disjoint union} of these spaces.
Let $X$ be a space satisfying $\ctblga$. Then the space $X\sqcup X$ satisfies $\ctblga$.
In the realm of sets of reals, the properties $\ctblga$ and $\GNga$ are equivalent.

\blem [{\cite[Proposition~2.3.]{miller}}]\label{m}
If $X, Y$ are sets of reals then the space $X\x Y$ satisfies  $\GNga$ if and only if the space $X\sqcup Y$ satisfies  $\GNga$.
\elem

From our main result we can obtain the following corollary which has originally was proved by Szewczak and W\l udecka~\cite[Theorem~4.1.(1)]{unbddtower}.

\bcor
Let $n\in\bbN$ and $X_1\cup \Fin,\dotsc, X_n\cup\Fin$ be  $\fp$-generalized tower sets  with the Cantor topology. 
Then the space $(X_1\cup\Fin)\x\dotsb\x(X_n\cup\Fin)$ satisfies $\ctblga$.
\ecor

\bpf
We prove the statement for $n=2$. The proof for other $n$ is similar.
Let $X, Y$ be $\fp$-generalized towers in $\roth$. 
Then $X\cup Y$ is a $\fp$-generalized tower.
By Theorem~\ref{thm1}, the space $X\cup Y\cup \Fin$ with the Michael line topology satisfies $\ctblga$.
Then the space $(X\cup Y\cup \Fin)\sqcup (X\cup Y\cup \Fin)$ satisfies $\ctblga$.
Since the property $\ctblga$ is hereditary for closed subset, thus the space $(X\cup\Fin)\sqcup (Y\cup\Fin)$ with the Michael line topology satisfies $\ctblga$.
Then $(X\cup\Fin)\sqcup (Y\cup\Fin)$ with the Cantor topology satisfies $\ctblga$, and $\GNga$, too. By Lemma~\ref{m}, the space $(X\cup \Fin)\x (Y\cup \Fin)$ with the Cantor topology satisfies $\GNga$.
\epf


\begin{thebibliography}{99}

\Pa{arch}{A. Arhangel'ski\u{\i}}{The frequency spectrum of a topological space and the classification of spaces}{Soviet Math. Dokl.}{13}{1972}{1185}{1189}

\Pa{arch2}{A. Arhangel'ski\u{\i}}{Hurewicz spaces, analytic sets and fan tightness of function spaces}{Soviet Mathematics Doklady}{33}{1986}{396}{399}


\bibitem{BarJu}
T.~Bartoszy\'nski, H.~Judah,
Set Theory: On the structure of the real line, A. K. Peters,
Massachusetts: 1995.

\Pa{BaTs}{T. Bartoszy\'nski, B. Tsaban}{Hereditary topological diagonalizations and the Menger--Hurewicz Conjectures}{Proceedings of the American Mathematical Society}{134}{2006}{605}{615}
	
\bibitem{blass} A. Blass, \emph{Combinatorial cardinal characteristics of the continuum}, in: \textbf{Handbook of Set Theory} (M. Foreman, A. Kanamori, eds.), Springer, 2010, 395--489.

\Pa{QN}{L. Bukovsk\'y, J. Hale\v{s}}{QN-spaces, wQN-spaces and covering properties}{Topology and its Applications}{154}{2007}{848}{858}

\Pa{wQN}{L. Bukovsk\'y}{On $\mathrm{wQN}_*$ and $\mathrm{wQN}^*$ spaces}{Topology and its Applications}{156}{2008}{24}{27}

\Pa{brp}{L. Bukovsk\'{y}, I.  Rec\l{}aw, M. Repick\'y}{Spaces not distinguishing convergences of real-valued functions}{Topology and its Applications}{112}{2001}{13}{40}

\Pa{gn}{J. Gerlits, Zs. Nagy}{Some properties of $\Cp(X)$, I}{Topology and its Applications}{14}{1982}{151}{161}
	
\Pa{gami}{F. Galvin, A. Miller}{$\gamma$-sets and other singular sets of real numbers}{Topology and its Applications}{17}{1984}{145}{155}

\Pa{hales}{J. Hale\v{s}}{On Scheepers' conjecture}{Acta Universitatis Carolinae. Mathematica et Physica}{46}{2005}{27}{31}

\Pa{coc2}{W. Just, A. Miller, M. Scheepers, P. Szeptycki}{The combinatorics of open covers II}{
	Topology and its Applications}{73}{1996}{241}{266}

\Pa{Laver}{R. Laver}{On the consistency of Borel's conjecture}{Acta Mathematicae}{137}{1976}{151}{169}

\bibitem{miller} A.~Miller, \emph{A hodgepodge of sets of reals}, Note di Matematica \textbf{27} (2007), suppl. 1, 25--39.

\Pa{BBC}{A. Miller, B. Tsaban}{Point-cofinite covers in Laver's model}{Proceedings of the American Mathematical Society}{138}{2010}{3313}{3321}

\Pa{gamma}{A. Miller, B. Tsaban, L. Zdomskyy}{Selective covering properties of product spaces, II: $\gamma$ spaces}{Transactions of the American Mathematical Society}{368}{2016}{2865}{2889}

\Pa{ot}{T. Orenshtein, B. Tsaban}{Linear $\sigma$-additivity and some applications}{Transactions of the American Mathematical Society}{363}{2011}{3621}{3637}	

\bibitem{sss}A. Osipov, P. Szewczak, B. Tsaban, \emph{Strongly sequentially separable function spaces, via selection principles}, Topology and its Applications, \textbf{270} (2020), 106942.

\bibitem{sakai} M.~Sakai, \emph{Property C'' and function spaces}, Proceedings of the American Mathematical Society \textbf{104} (1988), 917--919. 

\Pa{sakaiCp}{M. Sakai}{The sequence selection properties of $\Cp(X)$}{Topology and its Applications}{154}{2007}{552}{560}

\Pa{sakaisemcont}{M. Sakai}{Selection principles and upper semicontinuous functions}{Colloquium Mathematicum}{117}{2009}{251}{256}.

\bibitem{SakaiScheepersPIT} M. Sakai, M. Scheepers,
\emph{The combinatorics of open covers}, in:
\textbf{Recent Progress in General Topology III}
(K. Hart, J. van Mill, P. Simon, eds.),
Atlantis Press, 2014, 751--800.

\Pa{coc1}{M. Scheepers}{Combinatorics of open covers. I: Ramsey theory}{Topology and its Applications}{69}{1996}{31}{62}

\Pa{SchCp}{M. Scheepers}{Sequential convergence in $\Cp(X)$ and a covering property}{East-West Journal of Mathematics}{1}{1999}{207}{214}

\Pa{alpha_i}{M. Scheepers}{$\Cp(X)$ and Arhangel'ski\u{\i}'s $\alpha_i$ spaces}{Topology and its Applications}{89}{1998}{265}{275}

\Pa{CBC}{M. Scheepers, B. Tsaban}{The combinatorics of Borel covers}{Topology and its Applications}{121}{2002}{357}{382}

\bibitem{ST} P.~Szewczak, B.~Tsaban, \emph{Products of Menger spaces: A combinatorial approach}, Annals of Pure and Applied Logic \textbf{168} (2017), 1--18.

\Pa{unbddtower}{P. Szewczak, M. Włudecka}{Unbounded towers and products}{Annals of Pure and Applied Logic}{172}{2021}{102900}


\bibitem{add} B.~Tsaban, \emph{Additivity numbers of covering properties}, in: \textbf{Selection Principles and Covering Properties in Topology} (L.~Kocinac, editor), Quaderni di Matematica 18, Seconda Universita di Napoli, Caserta 2006, 245--282.

\bibitem{MHP} B. Tsaban,
\emph{Menger's and Hurewicz's Problems: Solutions from ``The Book'' and refinements},
Contemporary Mathematics \textbf{533} (2011), 211--226.

\Pa{sfh}{B. Tsaban, L. Zdomskyy}{Scales, fields, and a problem of Hurewicz}{Journal of the European Mathematical Society}{10}{2008}{837}{866}

\Pa{TsZdArh}{B. Tsaban, L. Zdomskyy}{Hereditarily Hurewicz spaces and Arhangel'ski\u{\i} sheaf amalgamations}{Journal of the European Mathematical Society}{12}{2012}{353}{372}
	
\end{thebibliography}
\end{document}